\documentclass[12pt]{report}

\usepackage{amsmath,amsfonts,graphicx,color,hyperref}

\usepackage{geometry}
\usepackage{fullpage}
\usepackage{epsfig}
\renewcommand{\baselinestretch}{1.3} %1.7 makes doublesp. wants 3 lines/2.5cm
\newcommand{\qed}{\nobreak \ifvmode \relax \else
      \ifdim\lastskip<1.5em \hskip-\lastskip
      \hskip1.5em plus0em minus0.5em \fi \nobreak
      \vrule height0.75em width0.5em depth0.25em\fi}
\begin{document}
\def\hk{{\mathcal H}_K}
%%to change   \cV to {\mathcal V}, use  :%s/\\\cV/{\\\mathcal V}/g
%%to change \tt to {\mbf{t}}, use :%s/\\\tt/{\\\mbf{t}}/g 
\def\cA{{\mathcal A}}
\def\cB{{\mathcal B}}
\def\cC{{\mathcal C}}
\def\cE{{\mathcal E}}
\def\cF{{\mathcal F}}
\def\cH{{\mathcal H}}
\def\cI{{\mathcal I}} 
\def\cL{{\mathcal L}}
\def\cM{{\mathcal M}}
\def\cM{{\mathcal M}}
\def\cN{{\mathcal N}}
\def\cS{{\mathcal S}}
\def\cT{{\mathcal T}}
\def\cU{{\mathcal U}}
\def\cX{{\mathcal X}}
\def\cY{{\mathcal Y}}
\def\ss{\mbf{s}}
\def\hk{{\mathcal H}_K}
\def\al{\alpha}
\def\be{\beta}
\def\ga{\gamma}
\def\de{\delta}
\def\ep{\epsilon}
\def\ze{\zeta}
\def\et{\eta}
\def\th{\theta}
\def\ka{\kappa}
\def\la{\lambda}
\def\rh{\rho}
\def\si{\sigma}
\def\ta{\tau}
\def\up{\upsilon}
\def\ph{\phi}
\def\ch{\chi}
\def\ps{\psi}
\def\om{\omega}
\def\epv{\vareplison}
\def\thv{\vartheta}
\def\Ga{\Gamma}
\def\De{\Delta}
\def\Th{\Theta}
\def\La{\Lambda}
\def\Si{\Sigma}
\def\tSi{\tilde{\Sigma}}
\def\tS{\tilde{S}}
\def\tM{\tilde{M}}
\def\tOm{\tilde{\Omega}}
\def\Up{\Upsilon}
\def\Ph{\Phi}
\def\Ps{\Psi}
\def\Om{\Omega}
\def\pa{\partial}
\def\del{\partial}
\def\beq{\begin{equation}}
\def\eeq{\end{equation}}
\def\mpar{\marginpar}
\def\epar{\marginpar} %use for equations 
\def\lrarrow{\leftrightarrow}
\renewcommand{\baselinestretch}{1.3} %1.7 makes doublesp. wants 3 lines/2.5cm
%\newcommand{\qed}{\nobreak \ifvmode \relax \else
      %\ifdim\lastskip<1.5em \hskip-\lastskip
      %\hskip1.5em plus0em minus0.5em \fi \nobreak
      %\vrule height0.75em width0.5em depth0.25em\fi}
%%%%%%%%%%%%%%%%%
\newcommand{\mbf}[1]{\mbox{\boldmath $#1$}}
\small \vspace{-.3in}
\thispagestyle{empty}
\setcounter{page}{0}
\textcolor{blue}
  {\flushleft{DEPARTMENT OF STATISTICS}\\
  University of Wisconsin\\
  1300 University Ave.\\
  Madison, WI 53706}
%\vfill %this was out
\vspace{1.3in}
%\vspace{.3in}
\begin{center}\large
    {\large TECHNICAL REPORT NO. 1173}\\
%\bigskip          % report number
     \today \\
% \vfill
\vspace{.1in}
\end{center}
\begin{center}
\Large{\bf{\sf \textcolor{blue}{
Statistical Model Building, Machine Learning, and the Ah-Ha Moment}}}
\end{center}
\vspace{.1in}
\begin{center}
{\sf \textcolor{blue}{Grace Wahba}\footnote[1]{Research supported in part by
NIH Grant EY09946 and NSF Grant DMS-0906818.
Prepared for the volume "Past, Present and Future of Statistics, 
in Celebration of the 50th Anniversary of the Committee of 
Presidents of Statistical Societies (COPSS)
}\\
Department of Statistics, Department of Computer Sciences \\ and 
Department of Biostatistics and Medical Informatics\\
University of Wisconsin, Madison
\vspace{.05in}
}
\end{center}
\normalsize
\newpage
\begin{center}
Preamble
\end{center}

As a recipient of the Elizabeth Scott Award, I 
received the following invitation from David Scott and COPSS:

``The Committee of Presidents of Statistical Societies (COPSS) 
will celebrate its 50th 
Anniversary in 2013. As part of its celebration, 
COPSS intends to publish a book with 
contributions from the past recipients of its four awards, 
namely the Fisher Lecture 
Award, the President's Award, the Elizabeth Scott Award, 
and the FN David Award. The 
theme of the book is “Past, Present and Future of Statistical Science.”
As a past winner of one of the COPSS awards, we would like 
to invite you to contribute 
to this book...
I will be working with you personally on your contribution ...
Specifically, we are seeking contributions along (one of) the following lines:
\begin{enumerate}
\item Statistical Research:  Your  reflection  on the past, 
present or future of a 
statistical research area of your choice. You are free to 
decide what you would 
like to focus on, i.e., past, present or future, or possibly all three.

\item  Statistical Education: Your
reflection/view on  statistical education in some 
areas, e.g., BIG DATA, interdisciplinary research,  graduate and undergraduate 
curriculum, promoting statistical education and 
research in developing countries, 
or statistical educational outreach.

\item  Statistical Career: Your reflection on your own career, 
lessons and experience 
you have learned, and advice you would like to 
provide to young statisticians if 
sought.

\item A blend of the above three topics. "
\end{enumerate}
(signed by David Scott representing COPSS). 
\\
\\
I have chosen to focus on Item 3,  reflecting on 
some particular ``fun" 
pieces of my own statistical career,
sprinkled with advice, 
and augmented with a few assorted remarks.
Following is my contribution.
\chapter{Statistical Model Building, Machine Learning, and the 
Ah-Ha Moment}
[\today ~ Chapter by Grace Wahba for the Committee of Presidents of 
Statistical Societies (COPSS) 50th anniversary volume]

Highly selected ``Ah-Ha" moments from the beginning to
the present  of my research career are 
recalled - these are moments when the main idea just 
popped up instantaneously, sparking sequences of future
research activity- 
almost all of these moments crucially
involved discussions/interactions  with others.
Along with a description of these moments we 
give unsought advice to young statisticians.
We conclude  
with remarks on issues relating to 
statistical model building/machine 
learning in the context of human subjects data. 

\section{Introduction-Manny Parzen and RKHS}\label{intro}
Many of  the ``Ah-Ha" moments below involve Reproducing
Kernel Hilbert Spaces (RKHS) so we begin there. 
My introduction to RKHS  
came while attending a class given by Manny Parzen 
on the lawn in front of the old Sequoia Hall at 
Stanford around 1963. See \cite{parzen:1962}.

For many years RKHS 
\cite{aronszajn:1950, wahba:1990}
were a little niche
corner of research which suddenly became popular when 
their relation to Support Vector Machines (SVMs)  became clear-
more on that later. To understand most of the Ah-Ha moments 
it may help to know a few facts about RKHS which
we now give.

An RKHS is a Hilbert space $\cH$ where all of the 
evaluation functionals are bounded linear functionals. 
What this means is the following: Let the domain of 
$\cH$ be $\cT$, and the inner product $<\cdot, \cdot>$.
Then, for  each $t \in \cT$ there exists an element, call it 
$K_t$ in $\cH$, with the property $f(t) = <f, K_t>$ for all 
$f$ in $\cH$. $K_t$ is known as the representer of 
evaluation at $t$. Let  $K(s,t) = <K_s, K_t>$; this is 
clearly a positive definite function on $\cT \otimes \cT$. 
By the Moore-Aronszajn theorem, every RKHS is associated
with a unique positive definite function, as we have
just seen. Conversely, 
given a positive definite function, there exists a unique
RKHS (which can be constructed from linear combinations
of the $K_t, t \in \cT$ and their limits). Given $K(s,t)$ we 
denote the associated RKHS as $\hk$. Observe
that nothing has been assumed concerning the domain 
$\cT$. A second role of positive definite functions 
is as the covariance of a zero mean Gaussian 
stochastic process on $\cT$. 
In a third role role that we 
will come across later - let  $O_i, i = 1,2,\dots, n$
be a set of $n$ abstract objects.
An $n \times n$ positive
definite matrix can be used to assign 
pairwise squared Euclidean distances 
$d_{ij}$ between $O_i$ and $O_j$ 
by 
$d_{ij}  = K(i,i) + K(j,j)- 2(K(i, j)$.
In Sections \ref{kimel}-\ref{dcor} we go 
through some Ah-Ha moments involving 
RKHS, positive definite functions and pairwise
distances/dissimilarities.  
Section \ref{reg} discusses sparse models and 
the LASSO. Section \ref{privacy} has some 
remarks involving complex interacting attributes,
the ``Nature-Nurture" debate, Personalized Medicine,
Human subjects privacy and scientific literacy,
and we end with conclusions in Section \ref{conclu}.

I end this section by noting that Manny Parzen was my 
thesis advisor, and Ingram Olkin was on my committee. 
My main advice to young statisticians is: Choose your 
advisor and committee carefully, and be as lucky as I 
was. 
%Another  member of my committee was Karel deLeeuw, 
%who sadly was murdered by a Stanford Math student who had
%spent 19 years unsuccessfully pursuing a PhD. 

\subsection{George Kimeldorf and the Representer Theorem}\label{kimel}
Back around 1970 George Kimeldorf and I both 
got to spend a lot of time at the Math Research 
Center at the University of Wisconsin-Madison
(the one that later got blown up as part
of the anti-Vietnam-war movement). At that time
it was a hothouse of spline work, headed by 
Iso Schoenberg, Carl deBoor, Larry Schumaker and others, and 
we thought that smoothing splines would be of 
interest to statisticians. The smoothing spline
of order $m$ was the solution to: find $f$ in the 
space of functions with 
square integral $m$th derivative to 
minimize

\begin{equation}\label{spline}
\sum_{i=1}^n (y_i - f(t_i))^2 + \lambda \int_0^1 (f^{(m)}(t))^2dt, 
~~~~t_i,~ i = 1, \dots ,n \in [0,1].
\end{equation}

Professor Schoenberg many years ago had characterized the 
solution to this problem as a piecewise polynomial 
of degree $2m-1$ satisfying some boundary and continuity
conditions. 

Our Ah-Ha moment came when we observed that 
the space of functions with square integrable $m$th 
derivative on $[0,1]$ was an RKHS 
with {\it seminorm} $\Vert Pf\Vert$ defined by   
$\Vert Pf\Vert^2 = \int_0^1 (f^{(m)}(t))^2 dt$ and with 
an associated $K(s,t)$ that we could figure out. 
(A seminorm is exactly like a norm except 
that it has a non-trivial null space, here the null space
of this seminorm is the span of the polynomials of 
degree $m-1$ or less.)
Then 
by replacing $f(t)$ by $<K_t, f>$ it was 
not hard to show by a very simple geometric
argument that the minimizer of (\ref{spline}) 
was in the 
span of the $K_t, t = t_1, \dots,  t_n$  and a 
basis for the  null space of 
the seminorm.  But furthermore, the very same
geometrical argument could be used to 
solve the more general problem: 
find $f \in \hk$, an  RKHS,  to minimize 
\begin{equation}\label{penlik}
\sum_{i=1}^n C(y_i, L_if) + \la\Vert Pf\Vert^2_K
\end{equation}
where $C(y_i,L_if)$ is convex in $L_if$, 
with $L_i$ a bounded linear functional in $\hk$
and $\Vert Pf \Vert^2_K$  a seminorm in $\hk$.
A bounded linear functional is a linear functional
with a representer in $\hk$, that is, there 
exists $\eta_{i} \in \hk$ such that 
$L_if = <\eta_{i}, f>$ for all $f\in \hk$. 
The minimizer of (\ref{penlik}) is in 
the span of the representers
$\eta_{i}$ and a basis for the null space
of the seminorm. 
That is known as the representer
theorem, which turned out to be a key to fitting 
(mostly continuous) functions in an infinite dimensional space, 
given a finite number of pieces of information.
There were two things I remember 
about our excitement over the result: One of us, I'm 
pretty sure it was George, thought the result was 
too trivial and not worthwhile to submit, but submit it we did and 
it was accepted
\cite{kimeldorf:wahba:1971} 
without a single complaint, within
three weeks. I have never since then had another paper accepted
by a refereed journal within three weeks and without
a single complaint. Advice: If you think it is worthwhile, 
submit it. 

\subsection{Svante Wold and Leaving-Out-One}
Following Kimeldorf and Wahba, it was clear that for practical
use, a method was needed to choose the smoothing or tuning 
parameter $\lambda$ in (\ref{spline}). The natural goal was to 
minimize the mean square error over the function $f$, for which 
its values at the data points would be the proxy.
In 1974 Svante Wold visited Madison, 
and we got to mulling over how to choose  $\lambda$. It so happened
that Mervyn Stone gave a colloquium talk in  Madison, and Svante and I were
sitting next to each other as Mervyn described using leaving-out-one
to decide on the degree of a polynomial to be used in least 
squares regression. We looked at each other at that very minute 
and simultaneously said something, I think it was  ``Ah-Ha", 
but possibly ``Eureka". In those 
days computer time was \$600/hour and Svante wrote a computer 
program to demonstrate that leaving-out-one did a good job. 
It took the entire Statistics department's computer money for an entire  month
to get the results in \cite{ wahba:wold:1975}. Advice: Go to the 
colloquia, sit next to your research pals. 

\subsection{Peter Craven, Gene Golub and Michael Heath and GCV}
After much struggle to prove some optimality properties
of leaving-out-one, it became clear that it couldn't 
be done in general. 
Considering the data model $y = f + \ep$, where
$y = (y_1, \dots, y_n)^T,  f = (f(t_1), \dots , f(t_n)^T$ 
and $\ep = (\ep_1, \dots, \ep_n)^T$ 
is a zero mean i.i.d. Gaussian random $n$-vector
then the information in the data is unchanged  
by multiplying left and right hand side 
by an orthogonal matrix, since $\Gamma\ep$ with $\Gamma$
orthogonal is still white Gaussian noise. 
But leaving-out-one can give you a different answer. 
To explain, we define the {\it influence matrix}: 
Let $f_{\lambda}$ be the minimizer of (\ref{spline}) when 
$C$ is sum of squares. The influence matrix
relates the data to the prediction of the data, 
$f_{\la} = A(\lambda)y$, where 
$f_{\la} = ( f_{\la}(t_1), \dots, f_{\la}(t_n))$.
A heuristic argument fell out of the blue, 
probably in an attempt to explain some
things to students, 
that rotating the data so that the influence 
matrix was constant down the diagonal, was the trick. 
The result was that instead of leaving-out-one, 
one should minimize  the GCV function 
$V(\lambda) = \frac{\sum_{i=1}^n(y_i -f(t_i))^2}{(trace(I-A(\lambda))^2}$
\cite{craven:wahba:1979,golub:heath:wahba:1979}.
I was on 
Sabbatical at Oxford in 1975 and Gene was at ETH 
visiting  Peter Huber, who had a beautiful house in 
Klosters, the fabled ski resort. Peter invited Gene and 
me up for the weekend, and Gene just wrote out the 
algorithm in \cite{golub:heath:wahba:1979} on the 
train from Zurich to Klosters while I snuck glances at
the spectacular scenery. Gene was a much loved mentor
to lots of people. He was born on February 29, 1932 and 
died in November of 2007. On February 29 and March 1, 2008 
his many friends held memorial birthday services at Stanford and 
30 other locations around the world. 
Ker-Chau Li \cite{li:1985, li:1986, li:1987b} and 
others later proved  optimality properties 
of the GCV and popular codes in R will compute splines
and other fits using GCV to estimate $\lambda$ and other
important tuning parameters. Advice: Pay attention to 
important tuning parameters since the results can be very 
sensitive to them. 
Advice: Appreciate mentors like Gene if you 
are lucky enough to have such great mentors.

\subsection{Didier Girard, Mike Hutchinson, 
Randomized Trace and the Degrees of Freedom 
for Signal}\label{df}
Brute force calculation of the trace of the influence matrix $A(\lambda)$ can 
be daunting to compute directly for large $n$. 
Let $f_{\lambda}^y$ be the minimizer of 
(\ref{spline}) with the data vector $y$ and let $f_{\la}^{y+\de}$ be 
the minimizer of (\ref{spline}) given the perturbed  data ${y + \de}$. 
Note that 
$\de^T (f_{\la}^y -f_{\la}^{y + \de}) = 
\de^TA(\lambda)(y + \de) -A(\lambda)(y)=
\sum_{i,j = 1}^n \de_i\de_ja_{ij}$, 
where $\de_i$ and $a_{ij}$ are the components of $\de$ and $A(\la)$ 
respectively. If the perturbations are i.i.d. with variance 1, 
then the expected value of this sum is an estimate of trace $A(\lambda)$.
This simple idea was proposed simultaneously
in \cite{girard:1989, hutchinson:1989}, with 
further theory in \cite{girard:1991}. It was a big Ah-Ha when I 
saw these papers because  further applications were immediate. In 
\cite{wahba:1983}, p 139, I defined the trace of $A(\la)$ as the 
``Equivalent degrees of freedom for signal", by analogy with 
linear least squares regression with  $p<n$  where the influence matrix is 
a rank $p$ projection operator. The degrees of freedom for signal
is an important concept in linear and nonlinear nonparametric 
regression, and it was a mistake to hide it inconspicuously in 
\cite{wahba:1983}. Brad Efron later \cite{efron:2004} gave an alternative
definition of degrees of freedom for signal.
The definition
in \cite{wahba:1983} depends only  on the data, Efron's is essentially
an expected value. Note that in the model (\ref{spline}), 
trace $A(\lambda) = 
 \sum_{i=1}^n \frac{\partial \hat{y_i}}{\partial y_i}$,
here $\hat{y_i}$ is the predicted value of $y_i$. 
This definition can reasonably be applied to a problem 
with a nonlinear forward operator (that is, that maps 
data onto the predicted data) when the derivatives 
exist, and the randomized trace method is reasonable 
for estimating the degrees of freedom for signal, 
although care should be
taken concerning the size of $\delta$. Even when 
the derivatives don't exist the randomized trace can 
be a reasonable way of getting at the degrees of 
freedom for signal, see for example 
\cite{wahba:johnson:gao:gong:1995}.  

\subsection{Yuedong Wang, 
Chong Gu and Smoothing Spline ANOVA}
Sometime in the late 80's or early 90's I 
heard Graham Wilkinson expound on 
ANOVA (Analysis of Variance), where 
data was given on a regular $d$-dimensional
grid: $y_{ijk},  t_{ijk}, i = 1, \dots, I, j = 1, \dots J, k = 1 \dots, K$, 
for $d = 3$ and so forth. That is, the domain is the  Cartesian product 
of several one-dimensional grids.  Graham was expounding on 
how fitting a model from observations on such a domain 
could be described as set of orthogonal projections
based on averaging operators, resulting in main effects, 
two factor interactions, etc. 
``Ah-Ha" I thought, we should be able to do exactly same thing 
and more where the domain  is the Cartesian product 
$\cT = \cT_1 \otimes \cT_2\otimes \dots \otimes \cT_d$
of  $d$  {\it arbitrary} domains. 
We want to fit functions on $\cT$, with main effects 
(functions of one variable), 
two factor interactions (functions of two variables), 
and possibly more terms given scattered observations, 
and we just need to define
averaging operators for each $\cT_{\alpha}$. 
Brainstorming with Yuedong Wang and Chong Gu fleshed
out the results. 
Let 
$\cH^{\al}, \al = 1, \dots, d$ be $d$ RKHSs  with 
domains $\cT_{\al}$, each $\cH^{\al}$  containing the constant
functions.
$\cH = \cH^1\otimes \cdots \otimes \cH^d$ is an RKHS with domain 
$\cT$. 
For each $\al=1,\dots, d$,  
construct a probability measure
$d\mu_{\al}$ on $\cT_{\alpha} $, with the property that 
the symbol  $(\cE_{\al}f)({t})$, the averaging operator, 
defined by 
\[
(\cE_{\al}f)({t})=\int_{{\cal{T}}^({\al})}f(t_1,\cdots,t_d)
d\mu_{\al}(t_\al),
\]
is well defined and finite for every $f \in \cH$ and ${t} \in 
\cT$. 
Consider the decomposition of the identity operator: 
\beq\label{E}
\noindent
I = \prod_{\al} (\cE_{\al} + (I-\cE_{\al})) 
\eeq
\beq
=\prod_{\al} \cE_{\al} + \sum_{\al}(I-\cE_{\al})\prod_{\be \neq \al} \cE_{\be}
+ \sum_{\al < \be} (I-\cE_{\al})(I-\cE_{\be})
\prod_{\ga \neq \al,\be}\cE_{\ga}
+\cdots+ \prod_{\al}(I - \cE_{\al}).
\eeq
This decomposition of the identity then always 
generates a unique (ANOVA-like)
decomposition of $f \in \cH$   of the form 
\begin{equation}\label{model}
f({t}) = \mu  + \sum_{\al} f_{\al}(t_{\al}) +
\sum_{\al < \be} f_{\al\be}(t_\al , t_\be ) 
+ \sum_{\al < \be < \ga}f_{\al\be\ga} (t_{\al}, t_{\be},t_{\ga}) +
\cdots
\end{equation}
where the expansion is  unique 
and (usually) truncated in some manner in practice. 
Here\\ 
$\mu  = (\prod_{\al}\cE_{\al})f, 
f_{\al} = ((I-\cE_{\al})\prod_{\be
\neq \al}\cE_{\be})f$, 
$f_{\al\be} = ((I -
\cE_{\al})(I-\cE_{\be})\prod_{\ga \neq \al,\be}\cE_{\ga})f $, etc, are
the mean, main effects, two factor interactions, etc. 
The result is usually called an SS ANOVA model, although 
the components are not limited to splines. 
For details on how to fit the terms see 
\cite{gu:2002,gu:wahba:1993,wahba:1990,wang:2011}
and the {\tt assist} and {\tt gss} codes in R. 
Note that {\it nothing} has been said about $\cT$ and 
very little regarding $\cH^{\al}$, other than that the 
constant functions are  in each of the constituent 
spaces and averaging 
operators can be defined. 

\subsection{Vladimir Vapnik, the Mystery Caller and the SVM}
%sve
The  AMS-IMS-SIAM Joint Summer Research Conference on
Adaptive Selection of Models and Statistical Procedures
was held on the campus
of Mount Holyoke College in South Hadley, Massachusetts on 
Sunday, June 23 
through Thursday, June 27, 1996.
On one of those fine days a session  met on a grassy
lawn of Mount Holyoke College, when Vladimir Vapnik
and I were both invited speakers. I talked first, 
and noted how the solution to the optimization 
problem (\ref{penlik}) led to a function involving the 
span of the representers. Vladimir spoke next, 
describing the support vector machine (SVM), a well
known and highly successful method for classification,
describing something he called the ``kernel trick".
He exhibited  an SVM that was fitted in 
the span of representers in an RKHS.
We will explain the SVM in a moment, but the 
original SVM, as proposed by Vapnik and 
coworkers \cite{vapnik:1995}  
was derived from an argument nothing like 
what I am about to give. Somewhere during 
Vladimir's talk, an unknown voice towards
the back of the audience called out 
``That looks like Grace Wahba's stuff." 
It looked obvious 
that the SVM as proposed by Vapnik with the 
``kernel trick", could be obtained as the
the solution to the optimization problem 
of (\ref{penlik}) with $C(y_i, L_i f)$ replaced 
by the so called hinge function, 
$(1-y_if(t_i))_+$, where $(\tau)_+ = \tau$ if $\tau> 0$ and 
$0$ otherwise. Each  data point  is  coded as $\pm 1$ 
according as it came from the ``plus" class or the 
``minus" class. For technical reasons the null space of 
the penalty function consists at most of the constant 
functions. 
Thus it follows that the solution is in the span of the representers 
$K_{t_i}$ from the chosen RKHS plus possibly a constant
function.  Yi Lin and coworkers
\cite{lin:lee:wahba:2002,lin:wahba:zhang:lee:2002} showed
that the SVM was {\it estimating the sign of the log 
odds ratio}, just what is needed for two class classification. 
The SVM may be compared to the case where one desires to 
estimate the {\it probability}  that an object is in the plus
class. If one begins with the penalized log likelihood 
of the Bernoulli distribution and codes the data 
as $\pm 1$ instead of the usual coding as $0$ or $1$, 
then we have the same optimization problem with 
$C(y_i, f(t_i)) = \log(1+e^{-y_if(t_i)})$
instead of $(1-y_if(t_i))_+$  
with solution in the same finite dimensional 
space, 
but it is estimating the log odds-ratio, 
as opposed to the {\it sign} of the log odds ratio. 
It was actually a big deal that the SVM could
be directly compared with 
penalized likelihood with Bernoulli data, 
and it provided a pathway for statisticians
and computer scientists to breach a major 
divide between them on the subject of 
classification, and to understand each others' work. 

For many years before the Hadley meeting, 
Olvi Mangasarian and I 
would talk about what we were doing in classification, 
neither of us having any understanding of what the 
other was doing. Olvi complained that the statisticians
dismissed his work, but it turned out that what he was
doing was related to the SVM and hence perfectly legitimate 
not to mention interesting, from
a classical statistical point of view. 
Statisticians and computer scientists have 
been on the same page on classification  ever since. 

It is curious to note that several patents have been 
awarded for the SVM. One of the early ones, issued 
on July 15, 1997 is 
``5649068 Pattern recognition system using support vectors"

I'm guessing that the unknown volunteer was David Donoho.

Advice: Keep your eyes open to synergies between 
apparently disparate fields.

\subsection{Yoonkyung Lee, Yi Lin and the Multicategory SVM}
%msvm.tex
For classification, when one has $k > 2$ classes 
it is always possible to apply an SVM to compare membership
in one class 
versus the rest 
of the $k$ classes, running through the algorithm $k$ times. 
In the early 2000s  there were many papers on one-vs-rest, 
and designs for subsets vs other subsets, but 
it is possible to generate examples where essentially
no observations will be identified as being in certain classes.
Since one-vs-rest could fail in certain circumstances
it was something of an open question how to do multicategory
SVMs in one optimization problem that did not have this 
problem. Yi Lin, Yoonkyung Lee 
and I were sitting around shooting the breeze and one of us 
said  ``how about a sum-to-zero 
constraint?" and the other two said ``Ah-Ha"!, or, 
at least that's the way I remember it. 
The idea is to code the labels as 
$k$-vectors, with a $1$ in the $r$th position and 
$-1/(k-1)$ in the $k-1$ other positions for a training sample
in class $r$. Thus, each observation 
vector satisfies the sum-to-zero constraint. 
The idea was to fit a vector of functions satisfying 
the same sum-to-zero constraint. 
The multicategory SVM fit estimates
$f(t) = (f_1(t), \cdots, f_k(t)), ~~ t \in \cT$ subject to 
the sum-to-zero constraint everywhere and the classification for 
a subject with attribute
vector $t$ is just the index of the largest component of 
the estimate of $f(t)$.
See \cite{lee:lee:2003,lee:lin:wahba:2004,lee:wahba:ackerman:2004}.
Advice: Shooting the breeze is good.

\subsection{Fan Lu, Steve Wright, Sunduz Keles,  Hector 
Corrada Bravo, and Dissimilarity Information}\label{rke}
We return to the alternative role of positive definite functions
as a way to encode pairwise distance observations. 
Suppose we are examining $n$ objects $O_i,~~ i = 1, \dots, n$ 
and are given some noisy or crude observations on their pairwise 
distances/dissimilarities, which may not satisfy
the triangle inequality. The goal is to embed these objects in 
a Euclidean space in such a way as to respect the 
pairwise dissimilarities as much as possible. Positive definite matrices
encode pairwise squared distances $d_{ij}$ between $O_i$ and $O_j$ as 
\begin{equation}
d_{ij}(K) = K(i,i) + K(j,j)- 2K(i,j),
\end{equation}
and, given a non-negative definite matrix of rank $d \leq n$, 
can be used to embed the $n$ objects in a Euclidean 
space of dimension $d$, centered at $0$ and unique up to rotations. 
We  seek  a $K$ which respects the 
dissimilarity information $d_{ij}^{obs}$ while 
constraining the complexity of $K$ by:
\begin{equation}\label{dijobs}
min_{K \in S_n} \sum |d_{ij}^{obs} -d_{ij}(K)| + \lambda~ trace K
\end{equation}
where $S_n$ is the convex cone of symmetric positive definite
matrices. 
I looked at this problem for an inordinate
amount of time seeking an analytic solution but after 
a conversation with Vishy (S. V. N.  Vishwanathan) at a meeting in 
Rotterdam in August of 2003 I realized it wasn't going 
to happen. The Ah-Ha moment came about when I showed the 
problem to Steve Wright, who right off said it could be solved
numerically using recently developed convex cone software. 
The result so far is 
\cite{corrada:lee:klein:klein:2009,lu:keles:wright:wahba:2005}.
In \cite{lu:keles:wright:wahba:2005} the objects are protein
sequences and the pairwise distances are BLAST scores. 
The fitted kernel $K$ had three eigenvalues that contained
about 95\% of the trace, so we reduced  $K$
to a rank 3 matrix by truncating the smaller eigenvalues. 
Clusters of four  different kinds of proteins were readily 
separated visually in three-d plots; see 
\cite{lu:keles:wright:wahba:2005} for the details.
In \cite{corrada:lee:klein:klein:2009} the objects are 
persons in pedigrees in a demographic study
and the distances are based on Malecot's kinship coefficient, 
which defines a pedigree dissimilarity measure. 
The resulting kernel became part of an SS ANOVA model 
with other attributes of persons, and the model
estimates a risk related to an eye disease.
Advice: Find computer scientist friends.

\subsection{Gabor Szekely, Maria Rizzo, Jing Kong and Distance
Correlation}\label{dcor}
The last Ah-Ha experience that we report is similar to that 
involving the randomized trace estimate of 
Section \ref{df}, that is, the Ah-Ha moment 
came about upon realizing that a particular recent 
result was very relevant to what we were doing. In this 
case Jing Kong brought to my attention the 
important paper of  Gabor Szekely and Maria Rizzo
\cite{szekely:rizzo:2009}. Briefly, this paper considers
the joint distribution of two random vectors, $X$ and $Y$, 
say, and 
provides a test, called distance correlation
that it factors so that the two
random vectors are independent. Starting with 
$n$ observations from the joint distribution, let $\{A_{ij}\}$ 
be the collection of double-centered pairwise distances among 
the $n \choose 2$ $X$ components, and similarly for $\{B_{ij}\}$.
The statistic, called distance
correlation,  is the analogue of the usual sample correlation
between the $A$'s and $B$'s.
The special property of the test is that it is justified 
for $X$ and $Y$ in Euclidean $p$ and $q$ space for arbitrary
$p$ and $q$ with no further distributional assumptions. 
In a demographic study involving pedigrees 
\cite{kong:klein:klein:lee:2012},
we observed that pairwise distance in death age between 
close relatives was less than that of unrelated age 
cohorts. 
A mortality risk score for four lifestyle factors 
and another score for a group of diseases was developed
via SS ANOVA modeling, 
and significant distance correlation was found between 
death ages, lifestyle factors and family relationships,
raising more questions than it answers regarding the 
"Nature-Nurture" debate (relative role of genetics
and other attributes). 

We take this opportunity to make a few important remarks about 
pairwise distances/dissimilarities, primarily how one 
measures them can be important, 
and getting the ``right" dissimilarity can be 90\% of the 
problem. We remark that family relationships in 
\cite{kong:klein:klein:lee:2012} were  based on 
a monotone function of Malecot's kinship coefficient
that was different from the monotone function 
in \cite{corrada:lee:klein:klein:2009}. Here it was
chosen to fit in with the different way the distances were used. In 
(\ref{dijobs}), the pairwise dissimilarities can be 
noisy, scattered, incomplete and could include subjective
distances like ``very close, close.. "  etc. not even satisfying
the triangle inequality. So there is substantial 
flexibility in choosing the dissimilarity measure with 
respect to the particular scientific context of the problem.
In \cite{kong:klein:klein:lee:2012} the pairwise distances
need to be a complete set, and be Euclidean (with some 
specific metric exceptions). There is still substantial
choice in choosing the definition of distance, since 
any linear transformation of a Euclidean coordinate
system defines a Euclidean distance measure.
Advice: Think about how you measure distance or 
dissimilarity in any problem involving pairwise relationships, 
it can be important. 

\section{Regularization Methods, RKHS and  Sparse Models}\label{reg}
The optimization problems in RKHS are a rich subclass of what can be 
called regularization methods, which solve an optimization 
problem which trades fit to the data 
versus complexity or constraints on
the solution.  My first encounter with the term 
``regularization" was \cite{tikhonov:1963} in the context
of finding  numerical  solutions to integral equations.
There the $L_i$ of (\ref{penlik}) were noisy integrals
of an unknown function one wishes to reconstruct, 
but the observations only contained a limited amount
of information regarding the unknown function. 
The basic and possibly revolutionary idea at the time was 
to find a solution which 
involves fit to the data while
constraining the solution by what amounted to 
an  RKHS seminorm, ($\int (f''(t))^2 dt$)
standing in for the missing 
information by an assumption that the solution was 
``smooth" \cite{o'sullivan:1986a,wahba:1977a}. 
Where once RKHS were  a niche subject, they are now
a major component of the statistical model building/machine 
learning literature. 

However, RKHS do not generally provide 
sparse models, that is,  models where a large number of coefficients
are being estimated but only a small but unknown 
number are believed to be  non-zero. 
Many problems in the ``Big Data" paradigm are believed to 
have, or want to have  sparse solutions, for example, 
genetic data vectors that may have many thousands of components 
and a modest number of subjects, as in a case-control 
study. 
The most popular method for ensuring sparsity is 
probably the LASSO 
\cite{chen:donoho:saunders:1998, tibshirani:1996}. 
Here  a very large 
dictionary of basis functions ($B_j(t),~~ j = 1, 2 \dots$)
is given and the unknown function is estimated as 
$f(t) =  \sum_j \beta_j B_j(t) $ with the penalty functional 
$\lambda\sum_{j} |\beta_j|$ replacing an RKHS square norm. 
This  will induce many zeroes
in the $\beta_j$, depending, among other things on the 
size of $\lambda$. Since then, researchers have commented
that there is a ``zoo" of proposed variants of sparsity-inducing 
penalties, many involving assumptions on  structures
in the data; one popular example is \cite{yuan:lin:2006}.
Other recent models involve mixtures of RKHS and sparsity-inducing
penalty functionals. 
One of our contribution to this ``zoo" deals with 
the situation where the data vectors amount  to very 
large ``bar codes", and it is desired to find 
patterns in the bar codes relevant to some outcome.
An innovative algorithm which deals with 
a humongous number of interacting patterns assuming that 
only a small number of coefficients are non-zero is given in 
\cite{shi:wahba:irizarry:corrada:2012,shi:wahba:wright:lee:2008,wright:2012}. 

As is easy to see here and in the statistical literature, 
the statistical modeler has overwhelming choices in modeling
tools, with many  public codes available in the software
repository R and elsewhere. In practice these choices
must be made with a serious understanding of the science 
and the issues motivating the data collection.
Good collaborations
with subject matter researchers can lead to the opportunity to 
participate in real contributions to the science. 
Advice: Learn
absolutely as much as you can about the subject matter of the 
data that you contemplate analyzing. When you use 
``black boxes" be sure you know what is inside them.

\section{Remarks on the Nature-Nurture Debate, 
Personalized Medicine and Scientific Literacy}\label{privacy}
We and many other researchers have been developing
methods for combining scattered, noisy, incomplete,
highly heterogenous information from multiple sources
with interacting variables to predict, classify, 
and determine patterns of attributes relevant to 
a response,  or more generally multiple 
correlated responses. 

Demographic studies, clinical trials, and ad hoc observational 
studies based on electronic medical records, which have 
familial \cite{corrada:lee:klein:klein:2009,kong:klein:klein:lee:2012},
clinical, genetic, lifestyle, treatment 
and other attributes
can be a rich source
of information regarding the Nature-Nurture 
debate, as well informing 
Personalized Medicine, two popular areas 
reflecting much activity. As large medical systems
put their records in electronic form interesting 
problems arise as to how to deal with such unstructured
data, to relate subject attributes to outcomes of interest.
No doubt a gold mine of information is 
there, particularly with respect to how the various attributes interact.
The statistical modeling/machine 
learning community continues to create and improve tools to 
deal with this data flood, eager to 
develop better and more efficient modeling methods, 
and regularization and dissimilarity methods will 
no doubt continue to play an important role in 
numerous areas of scientific endeavor.
With regard to 
human subjects studies, a limitation is  
the problem  of patient confidentiality - 
the more attribute information available to 
explore for its relevance, the trickier
the privacy issues, to the extent that 
de-identified data can actually be identified. 
It is important, however, that statisticians
be involved from the very start in the design of 
human subjects studies. 
%Procedures and methods for dealing with this difficult 
%situation which preserves privacy while 
%allowing researchers access to 
%data are very necessary.

With health related research, the US citizenry
has some appreciation of scientific results that 
can lead to better health outcomes. 
On the other hand any scientist who reads the 
newspapers or follows present day US politics 
is painfully aware that a non-trivial portion
of voters and the officials they 
elect 
have little or no understanding of 
the scientific method.
Statisticians need to participate 
in the promotion of 
increased scientific literacy 
in our educational establishment at all levels.

\section{Conclusion}\label{conclu}
In response to the invitation from COPSS to contribute
to their 50th Anniversary Celebration, I have taken a tour of 
some exciting moments in my career, 
involving RKHS and regularization methods, 
pairwise dissimilarities and distances, 
and LASSO models, dispensing un-asked for advice to new 
researchers along the way. I have made a few remarks concerning 
the richness of models based on RKHS, as well as models
involving sparsity-inducing penalties with some
remarks involving the Nature-Nurture Debate
and Personalized Medicine. 
I end this contribution with thanks to 
my many coauthors-identified here or not, and 
to my terrific present and former students. 
Advice: Treasure your collaborators! Have great students!

\bibliographystyle{plain}

\begin{thebibliography}{10}

\bibitem{aronszajn:1950}
N.~Aronszajn.
\newblock Theory of reproducing kernels.
\newblock {\em Trans. Am. Math. Soc.}, 68:337--404, 1950.

\bibitem{chen:donoho:saunders:1998}
S.~Chen, D.~Donoho, and M.~Saunders.
\newblock Atomic decomposition by basis pursuit.
\newblock {\em SIAM J. Sci. Comput.}, 20:33--61, 1998.

\bibitem{corrada:lee:klein:klein:2009}
H.~{{Corrada Bravo}}, K.~E. Lee, B.~E.~K. Klein, R.~Klein, S.~K. Iyengar, and
  G.~Wahba.
\newblock Examining the relative influence of familial, genetic and
  environmental covariate information in flexible risk models.
\newblock {\em Proceedings of the National Academy of Sciences},
  106:8128--8133, 2009.
\newblock Open Source at {\tt www.pnas.org/content/106/20/8128.full.pdf+html},
  PMCID: PMC 2677979.

\bibitem{craven:wahba:1979}
P.~Craven and G.~Wahba.
\newblock Smoothing noisy data with spline functions: estimating the correct
  degree of smoothing by the method of generalized cross-validation.
\newblock {\em Numer. Math.}, 31:377--403, 1979.

\bibitem{efron:2004}
B.~Efron.
\newblock The estimation of prediction error: covariance penalties and
  cross-validation.
\newblock {\em J. Amer. Statist. Assoc.}, 99:619--632, 2004.

\bibitem{girard:1989}
D.~Girard.
\newblock A fast `{M}onte-{C}arlo cross-validation' procedure for large least
  squares problems with noisy data.
\newblock {\em Numer. Math.}, 56:1--23, 1989.

\bibitem{girard:1991}
D.~Girard.
\newblock Asymptotic optimality of the fast randomized versions of ${GCV}$ and
  ${{C}_{{L}}}$ in ridge regression and regularization.
\newblock {\em Ann. Statist.}, 19:1950--1963, 1991.

\bibitem{golub:heath:wahba:1979}
G.~Golub, M.~Heath, and G.~Wahba.
\newblock Generalized cross validation as a method for choosing a good ridge
  parameter.
\newblock {\em Technometrics}, 21:215--224, 1979.

\bibitem{gu:2002}
C.~Gu.
\newblock {\em Smoothing Spline ANOVA Models}.
\newblock Springer, 2002.

\bibitem{gu:wahba:1993}
C.~Gu and G.~Wahba.
\newblock Smoothing spline {ANOVA} with component-wise {B}ayesian ``confidence
  intervals''.
\newblock {\em J. Computational and Graphical Statistics}, 2:97--117, 1993.

\bibitem{hutchinson:1989}
M.~Hutchinson.
\newblock A stochastic estimator for the trace of the influence matrix for
  {L}aplacian smoothing splines.
\newblock {\em Commun. Statist.-Simula.}, 18:1059--1076, 1989.

\bibitem{kimeldorf:wahba:1971}
G.~Kimeldorf and G.~Wahba.
\newblock Some results on {T}chebycheffian spline functions.
\newblock {\em J. Math. Anal. Applic.}, 33:82--95, 1971.

\bibitem{kong:klein:klein:lee:2012}
J.~Kong, B.~Klein, R.~Klein, K.~Lee, and G.~Wahba.
\newblock Using distance correlation and {S}moothing {S}pline {ANOVA} to assess
  associations of familial relationships, lifestyle factors, diseases and
  mortality.
\newblock {\em PNAS}, pages 20353--20357, 2012.

\bibitem{lee:lee:2003}
Y.~Lee and C.-K. Lee.
\newblock Classification of multiple cancer types by multicategory support
  vector machines using gene expression data.
\newblock {\em Bioinformatics}, 19:1132--1139, 2003.

\bibitem{lee:lin:wahba:2004}
Y.~Lee, Y.~Lin, and G.~Wahba.
\newblock Multicategory support vector machines, theory, and application to the
  classification of microarray data and satellite radiance data.
\newblock {\em J. Amer. Statist. Assoc.}, 99:67--81, 2004.

\bibitem{lee:wahba:ackerman:2004}
Y.~Lee, G.~Wahba, and S.~Ackerman.
\newblock Classification of satellite radiance data by multicategory support
  vector machines.
\newblock {\em J. Atmos. Ocean Tech.}, 21:159--169, 2004.

\bibitem{li:1985}
K.~C. Li.
\newblock From {S}tein's unbiased risk estimates to the method of generalized
  cross-validation.
\newblock {\em Ann. Statist.}, 13:1352--1377, 1985.

\bibitem{li:1986}
K.~C. Li.
\newblock Asymptotic optimality of ${{C}}_{{L}}$ and generalized cross
  validation in ridge regression with application to spline smoothing.
\newblock {\em Ann. Statist.}, 14:1101--1112, 1986.

\bibitem{li:1987b}
K.~C. Li.
\newblock Asymptotic optimality for {C} sub p , {C} sub {L} , cross-validation
  and generalized cross validation: discrete index set.
\newblock {\em Ann. Math. Statist.}, 15:958--975, 1987b.

\bibitem{lin:lee:wahba:2002}
Y.~Lin, Y.~Lee, and G.~Wahba.
\newblock Support vector machines for classification in nonstandard situations.
\newblock {\em Machine Learning}, 46:191--202, 2002.

\bibitem{lin:wahba:zhang:lee:2002}
Y.~Lin, G.~Wahba, H.~Zhang, and Y.~Lee.
\newblock Statistical properties and adaptive tuning of support vector
  machines.
\newblock {\em Machine Learning}, 48:115--136, 2002.

\bibitem{lu:keles:wright:wahba:2005}
F.~Lu, S.~Keles, S.~Wright, and G.~Wahba.
\newblock A framework for kernel regularization with application to protein
  clustering.
\newblock {\em Proceedings of the National Academy of Sciences},
  102:12332--12337, 2005.
\newblock Open Source at www.pnas.org/content/102/35/12332, PMCID: PMC118947.

\bibitem{o'sullivan:1986a}
F.~O'Sullivan.
\newblock A statistical perspective on ill-posed inverse problems.
\newblock {\em Statistical Science}, 1:502--527, 1986a.

\bibitem{parzen:1962}
E.~Parzen.
\newblock An approach to time series analysis.
\newblock {\em Ann. Math. Statist.}, 32:951--989, 1962.

\bibitem{shi:wahba:irizarry:corrada:2012}
W.~Shi, G.~Wahba, R.~Irizarry, H.~Corrada Bravo, and S.~Wright.
\newblock The partitioned {LASSO-P}atternsearch algorithm with application to
  gene expression data.
\newblock {\em BMC Bioinformatics}, 13-98, 2012.

\bibitem{shi:wahba:wright:lee:2008}
W.~Shi, G.~Wahba, S.~Wright, K.~Lee, B.~Klein, and R.~Klein.
\newblock {LASSO} {P}attern search algorithm with applications to ophthalmology
  and genomic data.
\newblock {\em Statistics and Its Interface}, 1:137--153, 2008.
\newblock SII-1-1-A12-Shi.pdf, PMCID:PMC2566544.

\bibitem{szekely:rizzo:2009}
G.~Szekely and M.~Rizzo.
\newblock Brownian distance covariance.
\newblock {\em Ann. Appl. Statist.}, 3:1236--1265, 2009.

\bibitem{tibshirani:1996}
R.~Tibshirani.
\newblock Regression shrinkage and selection via the {LASSO}.
\newblock {\em J. Roy. Stat. Soc, B}, 58:267--288, 1996.

\bibitem{tikhonov:1963}
A.~Tikhonov.
\newblock Solution of incorrectly formulated problems and the regularization
  method.
\newblock {\em Soviet Math. Dokl.}, 4:1035--1038, 1963.

\bibitem{vapnik:1995}
V.~Vapnik.
\newblock {\em The Nature of Statistical Learning Theory}.
\newblock Springer, 1995.

\bibitem{wahba:1977a}
G.~Wahba.
\newblock Practical approximate solutions to linear operator equations when the
  data are noisy.
\newblock {\em SIAM J. Numer. Anal.}, 14:651--667, 1977a.

\bibitem{wahba:1983}
G.~Wahba.
\newblock Bayesian ``confidence intervals'' for the cross-validated smoothing
  spline.
\newblock {\em J. Roy. Stat. Soc. Ser. B}, 45:133--150, 1983.

\bibitem{wahba:1990}
G.~Wahba.
\newblock {\em Spline Models for Observational Data}.
\newblock SIAM, 1990.
\newblock CBMS-NSF Regional Conference Series in Applied Mathematics, v. 59.

\bibitem{wahba:johnson:gao:gong:1995}
G.~Wahba, D.~Johnson, F.~Gao, and J.~Gong.
\newblock Adaptive tuning of numerical weather prediction models: randomized
  {GCV} in three and four dimensional data assimilation.
\newblock {\em Mon. Wea. Rev.}, 123:3358--3369, 1995.

\bibitem{wahba:wold:1975}
G.~Wahba and S.~Wold.
\newblock A completely automatic {F}rench curve.
\newblock {\em Commun. Stat.}, 4:1--17, 1975.

\bibitem{wang:2011}
Y.~Wang.
\newblock {\em Smoothing Splines: Methods and Applications}.
\newblock Chapman \& Hall/CRC Monographs on Statistics \& Applied Probability,
  2011.

\bibitem{wright:2012}
S.~Wright.
\newblock Accelerated block-coordinate relaxation for regularized optimization.
\newblock {\em S{IAM} J. Optimization}, 22:159--186, 2012.
\newblock Preprint and software available at {\tt
  http://pages.cs.wisc.edu/\~{}swright/LPS/}.

\bibitem{yuan:lin:2006}
M.~Yuan and Y.~Lin.
\newblock Model selection and estimation in regression with grouped variables.
\newblock {\em J. Roy. Statist. Soc. B}, 68:49--67, 2006.

\end{thebibliography}

\end{document}